\documentclass[12pt,twoside,a4paper]{amsart}
\usepackage{amssymb}
\date{\today}

\def\deg{\text{deg}\,}

\def\w{\wedge}

\def\dbar{\bar\partial}

\def\C{{\mathbb C}}
\def\P{{\mathbb P}}

\def\D{{\mathcal D}}

\def\codim{{\rm codim\,}}

\def\E{{\mathcal E}}

\def\O{{\mathcal O}}

\def\L{{\mathcal L}}
\def\Q{{\mathcal Q}}
\def\Re{{\rm Re\,  }}

\def\L{{\mathcal L}}

\def\be{\begin{equation}}
\def\ee{\end{equation}}

\newtheorem{thm}{Theorem}[section]
\newtheorem{lma}[thm]{Lemma}
\newtheorem{cor}[thm]{Corollary}
\newtheorem{prop}[thm]{Proposition}

\theoremstyle{definition}

\theoremstyle{remark}

\newtheorem{preremark}{Remark}
\newtheorem{preex}{Example}

\newenvironment{remark}{\begin{preremark}}{\end{preremark}}
\newenvironment{ex}{\begin{preex}}{\end{preex}}

\numberwithin{equation}{section}

\begin{document}

\title[The membership problem for polynomial ideals in  \dots]
{The membership problem for polynomial ideals in terms of 
residue currents}

%%%%{\center \Huge PRELIMINARY VERSION !!!}

\date{\today}

\author{Mats Andersson}

\address{Department of Mathematics\\Chalmers University of Technology and the University of G\"oteborg\\S-412 96 G\"OTEBORG\\SWEDEN}

\email{matsa@math.chalmers.se}



\thanks{The author was
  partially supported by the Swedish Natural Science
  Research Council}

\begin{abstract}
We find a relation between the vanishing of a globally defined residue current
on $\P^n$ and solution of the membership problem with control of the polynomial
degrees. Several classical results appear as special cases,
such as Max N\"other's theorem,   and we also obtain a generalization of that theorem.
There are also  
connections to effective versions of the Nullstellensatz.
We also provide explicit integral representations of the
solutions.
\end{abstract}

\maketitle

\section{Introduction}

Let $F_1, \ldots, F_m$ be polynomials in $\C^n$ and let $\Phi$ be a polynomial
that vanishes on the common zero set of the $F_j$. By Hilbert's Nullstellensatz,
for some power $\Phi^\nu$ of $\Phi$,  one can find polynomials
$Q_j$ such that
\begin{equation}\label{busa}
\sum_j F_j Q_j=\Phi^\nu.
\end{equation}
A lot of attention has been paid to find effective versions, i.e., control of $\nu$ and
the degrees of $Q_j$ in terms of the degrees of $F_j$.
%%%
 The breakthrough was in  \cite{Br} where Brownawell obtained
bounds on $\nu$ and $\deg Q_j$ not too far from  the best possible, using 
a combination of algebraic and analytic methods, cf., Remark~\ref{brownawell}
below.
%%that \eqref{estupp} holds with $M=(n-1)d^{\min(m,n)}-1+d$ (assuming $\deg F_j\le d$) 
%%and by applying theorem he obtained a solution
%%$Q$ with $\deg Q_j\le n\min(m,n)d^{\min(m,n)}+\min(m,n)d$.
Soon after that Koll\'ar \cite{Koll} obtained  by purely algebraic 
 methods the following optimal result.

\smallskip
\noindent{\bf Theorem (Koll\'ar).}
{\it Let $F_1, \ldots, F_n$ and $\Phi$ be  polynomials in $\C^n$
of degrees $d_j$, and $r$, respectively,  and assume that $\Phi$ vanishes on the
common zero set of $F_j$. Then
(if $d_j\neq 2$), one can find 
polynomials $Q_j$ and a natural number $s$ such that
$\sum F_jQ_j=\Phi^\nu$, 
and such that
$\nu\le N(d_1\cdots d_m)$ and $\deg F_jQ_j\le (1+r)N(d_1\cdots d_m)$;
here
$N(d_1\cdots d_m)=d_1\cdots d_m$ if $m\le n$;
for the case when  $m>n$, see \cite{Koll}.}
\bigskip

In particular, if $F_j$ have no common zeros in $\C^n$, then there are
polynomials $Q_j$ such that
\begin{equation}\label{bezout}
\sum_j F_jQ_j=1,
\end{equation}
with
$$
\deg   F_jQ_j\le  N(d_1\cdots d_m).
$$
The restriction $d_j\neq 2$ has recently been removed by Jalonek, \cite{Jal},
in the case when $m=n$.
\smallskip

In \cite{Br2} Brownawell gave a prime power version of Koll\'ar's theorem
which shed more geometric light on these questions, 
and there is a  generalization  to smooth algebraic manifolds in \cite{EL}.
%%See also Sombra, ???.\, ???? etc for alternative proofs and related results.

\smallskip
Koll\'ar's  result is optimal  as long as one  only makes   assumptions of the degrees of $F_j$.
However,  if one imposes geometric conditions   on the zero set
one can get  sharper results.
For instance, assuming that $m=n+1$ and $F_j$ have no common zero set even at
infinity, then a classical theorem of Macaulay,  \cite{Macaul}, 
states that \eqref{bezout}  has a solution such that
$\deg F_jQ_j\le\sum d_j-n$.

\smallskip

There is a  related  result 
due to Max N\"other, \cite{Noe};  see also \cite{GH}.

\smallskip

\noindent{\bf Theorem  (Max N\"other, 1873).}
{\it Assume that the zero set of $F_1,\ldots, F_n$ is discrete and contained in $\C^n$
and that $\Phi$ belongs to the ideal $(F)$. Then there are polynomials $Q_j$ such that
$$ %% \begin{equation}\label{plums}
\Phi=\sum_1^n F_jQ_j 
$$ %% \end{equation}
and 
$\deg F_jQ_j\le \deg \Phi.
$
}

\bigskip

In this paper we present a more general result about solutions 
to the equation   
\begin{equation}\label{plums}
\Phi=\sum_1^m  F_jQ_j,
\end{equation}
where  $F_1,\ldots,F_m$ are given  polynomials in $\C^n$,
with control of the degrees of $F_jQ_j$.
It is formulated in terms of  a residue current associated with
$F_j$ with support on their common zero set on $\P^n$,
and  the theorems of Macaulay and Max N\"other are  simple 
consequences. 
We also provide explicit representation formulas of  solutions. 
%%%

\smallskip
If  $f_j$ denote  homogenizations of $F_j$, i.e., $f_j(z)=z_0^{d_j}F_j(z'/z_0)$,
where $d_j\ge\deg F_j$, (here $z=(z_0,z_1,\ldots,z_n)$ and $z'=(z_1,\ldots,z_n)$),
then each $f_j$ defines a global holomorphic section of the 
line bundle $L^{d_j}\to\P^n$, and hence
$f=(f_1,\cdots, f_m)$ is a section of 
the rank $m$ bundle $E^*=L^{d_1}\oplus\cdots\oplus L^{d_m}$ over $\P^n$
(here  $L^s$ denotes  the  line bundle $\O(s)$).
If $z\in\C^{n+1}\setminus\{0\}$ we let $[z]$ denote the corresponding
point in $\P^n$ under the natural projection; %%sometimes, 
however, we write  $f(z)$ rather than $f([z])$. 
%%Moreover,
%%$z'=(z_1,\ldots,z_n)$.)
%%
%%
If $E^*$ is equipped with the natural Hermitian structure,
then
\begin{equation}\label{pnorm}
\|f(z)\|^2=\sum_1^m\frac{|f_j(z)|^2}{|z|^{2d_j}}.
\end{equation}
Following \cite{A2} we  can define the  residue current $R^f$ 
which is an element in
$\oplus_{\ell}\D'_{0,\ell}(\P^n,\Lambda^\ell E)$ and  with support on
the zero set 
$$
Z^f=\{[z]\in\P^n;\ f(z)=0\}.
$$
If we assume that
the polynomials $F_j$ have no common zeros in $\C^n$, then of course $Z^f$ is
a subset of the hyperplane at infinity.
If $\codim Z^f=m$, i.e., $f$ is locally a complete intersection, then
$R^f$ is a  $(0,m)$-current
with values in $\det E=L^{-\sum d_j}$; more precisely it can be identified with the
 Coleff-Herrera current  %%%which formally can be written
$$
\Big[\dbar\frac{1}{f_1}\wedge\ldots\wedge\dbar\frac{1}{f_m}\Big],
$$
in $\C^{n+1}\setminus\{0\}$, see Section~\ref{rest}.
We can  now formulate our main result in this paper.

\begin{thm}\label{noll}
Let $F_1,\ldots, F_m$ be polynomials in $\C^n$, $\deg F_j\le d_j$,  
let $f=(f_1,\cdots,f_m)$ be the corresponding
section of   $E^*=L^{d_1}\oplus\cdots\oplus L^{d_m}$ over $\P^n$,
and let $R^f$ be the asssociated  residue current.
Moreover, assume that
\begin{equation}\label{valle}
\quad m\le n  \quad {\rm or}\quad r\ge \sum_{j=1}^{n+1} d_j-n,
\end{equation}
where $d_1\ge d_2\ge\ldots\ge d_m$. 
Let $\Phi$ be a polynomial, $\deg\Phi\le r$,  and let
$\phi\in\O(\P^n,L^r)$  denote its  $r$-homogenization.
If  
\begin{equation}\label{resvillkor}
\phi R^f=0,
\end{equation}
then there are polynomials $Q_j$ such that
\eqref{plums} holds 
%%\begin{equation}\label{snorvalp}
%%\Phi=\sum_1^m  F_jQ_j
%%\end{equation}
and  $\deg  F_jQ_j\le r$.
If $f$ is a complete intersection (then the condition \eqref{valle} is fulfilled) and
there exist   such polynomials  $Q_j$, then 
\eqref{resvillkor} holds.
\end{thm}

It is clear that the conclusion about $\deg F_jQ_j$ cannot be  improved.
%%%
If $\Phi=1$ the condition \eqref{resvillkor} means  that $F_j$ have no common zeros
in $\C^n$ and that $z_0^r$ annihilates  the residue
$R^f$ at infinity. 
%%%
If $Z^f$ is empty and $m=n+1$ (actually any $m\ge n+1$ works)  we can choose
$r=\sum d_j -n$ and hence we get a solution to the B\'ezout equation  \eqref{bezout}
such that $\deg  F_jQ_j\le \sum d_j-n$; thus we have obtained the
theorem of Macaulay mentioned above.

\smallskip
We  have  the following  generalization of
N\"other's theorem.

\begin{thm}\label{gennoth}
Assume that the projective zero set of $F_1,\ldots, F_m$  has codimension $m$ and that
there is no irreducible component contained in the hyperplane at infinity. If $\Phi$ belongs to
ideal $(F)$, then there are polynomials $Q_j$ such that
\eqref{plums} holds and
%%\begin{equation}%\label{plums}
%%\Phi=\sum_1^n F_jQ_j 
%%\end{equation}
$
\deg F_jQ_j\le \deg \Phi.
$
\end{thm}

\begin{proof}%%%%%%%%[Proof of Max N\"other's theorem and Theorem~\ref{gennoth}]
Since $m\le n$ the condition
\eqref{valle} is fulfilled so  we can take $r=\deg \Phi$. 
Since $\Phi\in(F)$,  $\phi$ is in the ideal $(f)$ locally in $\C^n$ and
since $f$ is a complete intersection, $\phi R^f=0$ in $\C^n$ by the 
duality theorem (see Section~\ref{rest}).
If $m=n$, i.e., as in  N\"other's theorem,   $Z^f$ is  contained in $\C^n\subset\P^n$,
so $R^f$ has its support in $\C^n$ as well, and hence
\eqref{resvillkor} holds in $\P^n$.
Thus Theorem~\ref{noll}   provides the desired solution.
In the general case the assumption means  that  the intersection of
$Z^f$ and the hyperplane at infinity has codimension $m+1$, and then 
Proposition~\ref{gubb}  in Section~\ref{rest} implies that
$\phi R^f=0$ in $\P^n$.
\end{proof}

\begin{remark}
Although  this theorem is probably known before, we have not found it in the literature.
A proof of N\"other's theorem by multivariable residue calculus has previously been obtained
by Tsikh, \cite{Ts}. In  \cite{TY} is given an argument starting with a representation of
$\Phi$ with the Cauchy-Weil formula.  Making series expansion of the kernel and using Jacobi
formulas (vanishing of certain residues as in  \cite{VY})  and the duality theorem,
one obtains N\"other's theorem. It is possible
 that one can prove  the  general form
of Theorem~\ref{gennoth} in a similar way, following the idea of
\cite{BT} oto  add $n-k$ linear forms $L$
such that $(F,L)$ has no zeros at infinity, but we have not checked the details.  
\end{remark}

However some results related to Theorem~\ref{gennoth}
have appeared before. In  \cite{Sh}, Proposition~2,
it is assumed   that  
$f_j$ is a regular sequence in $\P^n$ but with no extra condition on 
the hyperplane at infinity. If $\Phi$ belongs to the ideal $(F)$ as above, then
there are $Q_j$ solving \eqref{plums}  such that $\deg F_j\Psi_j\le N+\deg\Phi$, 
where $N=\Pi_1^m\deg F_j$.  To see this in our setting, 
recall  that (see, e.g., \cite{Sh} Lemma~2)  if $F_j$ is a regular sequence in $\O_x$, then
\begin{equation}\label{ink}
(\sqrt {(F)_x})^N\subset (F)_x.
\end{equation}
Thus $z_0^{N+\deg\Phi}\Phi(z'/z_0)$
annihilates $R^f$ in $\P^n$, and therefore the statement
 follows from Theorem~\ref{noll}.

\smallskip

Now let   $F_j$ be  as in Theorem~\ref{gennoth} and assume that
$\Phi$ vanishes on their common zero set  in  $\C^n$. Then by \eqref{ink},
  $\Phi^N$ belongs to
$(F)$. Therefore we get the following corollary of Theorem~\ref{gennoth}, which
recently appeared in \cite{FPT} under the slightly stronger assumption that
$F_j$ is a strictly regular sequence in $\C^n$.

\begin{cor}
Assume that the projective zero set of $F_1,\ldots, F_m$  has codimension $m$ and that
there is no component contained in the hyperplane at infinity. If $\Phi$ vanishes on
the zero set of $F$ in $\C^n$, then there are polynomials $Q_j$ such that
$
\deg F_jQ_j\le N\deg \Phi
$
and $\sum F_j Q_j=\Phi^N$, where
$N=(\deg F_1)\cdots(\deg F_m)$.
\end{cor}

If 
\begin{equation}\label{strahatt}
\|\phi\|\le\ C\|f\|,
\end{equation}
then, see Section~\ref{rest},  $\phi^{\min(m,n)} R^f=0$,  and hence 
Theorem~\ref{noll} implies

\begin{cor}\label{skolk}
Let  $F_j$ and $\Phi$ be as in Theorem~\ref{noll}, $r\ge\deg\Phi$, 
and assume that
$$
m\le n \quad{\rm or}\quad r\min(m,n)\ge\sum_1^{n+1} d_j-n.
$$
If \eqref{strahatt} holds, then 
there are polynomials $Q_j$ such that
\begin{equation}\label{putta}
\sum F_jQ_j=\Phi^{\min(m,n)}
\end{equation}
and $\deg F_jQ_j \le r\min(m,n)$.
\end{cor}

Since there are examples where $f$ is a complete
intersection and the full power $\min(m,n)$ of $\phi$ is needed
to kill $R^f$,  this result is then sharp.

\begin{ex}\label{pantalon} 
Let  $M$ be a given positive integer. 
Take $F_j(z')=z_j^{Mm}$ in $\C^n$, $1\le j\le m\le n$,
 $\Phi(z')=(z_1+\cdots z_m)^{Mm}$, and let $f_j$ and $\phi$ be the
homogenizations as before ($d_j=Mm$).
Then $\eqref{strahatt}$ holds and hence the corollary states that \eqref{putta} has a 
solution such that $\deg F_jQ_j\le r\min(m,n)=Mm^2$.
This is obvious also by a direct inspection, and  one also immediately sees
that $\Phi^{m-1}$ is not in the ideal $(F)$. Thus the corollary is sharp.

It follows that  $\phi^{m-1}R^f\neq 0$, and since 
 $f$ is a complete intersection in $\P^n$,
in fact $Z^f$ is the $n-m$-plane $\{[z]\in\P^n;\ z_1=\cdots z_m=0\}$,
it also follows that  $\phi^mR^f=0$.
One can also verify these residue conditions directly. In fact,
in the standard affine coordinates $z'$,
$$
R^f=\Big[\dbar\frac{1}{z_1^{Mm}}\w\ldots\w\dbar\frac{1}{z_m^{Mm}}\Big]\w \epsilon,
$$
where $\epsilon$ is a non-vanishing section of the line bundle $\det E$,
see Section~\ref{rest}.
Since this residue current is a tensor product of one-variable currents,
the residue conditions  follow from the 
 one-variable equality
$z\dbar[1/z^{p+1}]=\dbar[1/z^p]$.
\end{ex}

\bigskip

Let $F_j$ be polynomials  with no common zeros in $\C^n$. 
Since the zero set of the section
$f$ (take $d_j=\deg F_j$) is then contained in the hyperplane at infinity
it follows from Lojasiewicz' inequality that
\begin{equation}\label{estuppx}
\|z_0\|^M\le C\|f\|
\end{equation}
for some  $M$, or equivalently,
%%holomorphic even across the hyperplane at infinity in $\P^n$ it follows  that
\begin{equation}\label{estupp}
\sum_1^m\frac{|F_j(z')|^2}{(1+|z'|^2)^{d_j}}\ge c\frac{1}{(1+|z'|^2)^M}.
\end{equation}
Under  this condition 
$z_0^{M \min(m,n)} R^f=0$,  so we have

\begin{cor}\label{oppo}
Let $F_1,\ldots, F_m$ be polynomials in $\C^n$ of degrees
$d_j$ such that \eqref{estupp} (or equivalently \eqref{estuppx})
holds for some number $M$, and  assume that 
$$
m\le n \quad{\rm or}\quad M\min(m,n)\ge\sum_1^{n+1} d_j-n.
$$
Then there is a solution to
$\sum F_jQ_j=1$ with
$\deg F_jQ_j\le \min(m,n)M$.
\end{cor}

\smallskip

\begin{ex}
Also Corollary~\ref{oppo} is essentially sharp.  Let $M$ be a given
non-negative integer and take  $F_j(z')=z_j^M$, 
$1\le j <m\le n$,  and $F_m(z')=(1+z_1+\cdots +z_{m-1})^M$.
Then $f_j=z_j^M$ and $f_m=(z_0+\cdots +z_{m-1})^M$, so
\eqref{estupp} holds. 
The corollary thus gives a solution to \eqref{bezout} with
$\deg F_jQ_j\le mM$.  

Writing $1=(1+z_1+\cdots +z_{m-1})-z_1-\cdots -z_{m-1}$ and taking the
power $Mm-m+1$ we get a solution to \eqref{bezout} with
$\deg F_jQ_j=Mm-m+1$, and it is easily seen to be the best possible,
cf., Example~\ref{pantalon}.  However,
for large $M$,  $Mm-m+1$ is close to
$Mm$.
\end{ex}

\begin{remark}
Given the estimate \eqref{estupp},
 one can obtain a solution to the B\'ezout equation \eqref{bezout}
by a direct application of Skoda's $L^2$-estimate from \cite{Sk}, as is done
in \cite{Br}. If we for simplicity assume that all $d_j=d$, then
one gets  a solution with $\deg Q_j\le \min(m+1,n)M-d$, i.e.,
$\deg F_jQ_j\le\min(m+1,n)M$.
For $m\le n$ this is  the same as in  Corollary~\ref{oppo} but for
$m>n$ it is strictly weaker. 
This phenomenon is the same as in the original proof
of  Brian\c con-Skoda's theorem, \cite{BS}. Under the assumption $\|\phi\|\le C\|f\|$,
the (local)   $L^2$-estimate  immediately implies that
$\phi^{\min(m,n+1)}$ belongs the  the ideal $(f)$ locally;
to obtain the correct result when $m>n$,  that the power
$n$ is enough,  an additional argument is required.
See also Section~\ref{rest}.
\end{remark}

\bigskip

\begin{remark}%%[Remark on Brownawell's result]
\label{brownawell}
The main step in  Brownawell's paper \cite{Br}
is to obtain good control of the power $M$ in
\eqref{estupp} in terms of the degrees
of $F_j$, assuming that they  have no common zeros
in $\C^n$, and this is done  by means of Chow forms,
see also \cite{Te}.

\smallskip
Koll\'ar's theorem implies  that the estimate \eqref{estupp} holds
with $M=N(d_1,\ldots,d_m)$, see, \cite{Koll},  and this is in fact best possible.
From this estimate one gets, via  Corollary~\ref{oppo}, a solution to \eqref{bezout}
with $\deg F_jQ_j\le \min(m,n)M$. In view of  Koll\'ar's theorem one has then
``lost'' the factor  $\min(m,n)$.
\end{remark}

\smallskip

\begin{remark}
Koll\'ar's theorem holds for any field.
Berenstein and Yger, \cite{BY}, have obtained
explicit solutions to the B\'ezout equation \eqref{bezout}
in subfields of
$\C$, by means of integral formulas; 
see also \cite{BGVY} and the more recent survey article
\cite{TY} for a thorough discussion.
\end{remark}

\smallskip

\begin{remark}
The condition \eqref{strahatt} means that $\phi$  locally on
$\P^n$ belongs to the integral closure of the ideal $(f)$.
In \cite{Hick}, Hickel proves that if $\Phi$ is in the integral closure
of $(F)$ in $\C^n$, then one can solve (assuming $m\le n$ for simplicity)
$\Phi^m=\sum F_jQ_j$ with $\deg(F_jQ_j)\le m\deg\Phi+m d_1\cdots d_m$.
This result would follow from Theorem~\ref{noll} if one could prove that
the current
$
z^{md_1\cdots d_m}_0\phi^m R^f
$
vanishes ($\phi$ is the $\deg\Phi$ homogenization of $\Phi$). 
In $\C^n$ it vanishes since $|\Phi|\le C|F|$ locally.
If the zero set is contained in $\{z_0=0\}$,  the current  vanishes
there by Koll\'ar's theorem. We do not know whether it vanishes 
in the general case.
%%how  one can see the
%%general case.
\end{remark}

Theorem~\ref{noll} is a special case of the following more general
result, for which we formulate only  the homogeneous
version. Let $\delta_f$ denote the
mapping $\E(\P^n,\Lambda^{\nu+1}E\otimes L^r)\to
\E(\P^n,\Lambda^{\nu}E\otimes L^r)$ defined as 
interior multiplication with the section $f$ of  $E^*$.
Thus for instance, if 
$q=(q_1,\ldots,q_m)$ is a section to $E\otimes L^r$, then
$\delta_f q$ is equal to the section $\sum_j f_jq_j$
of  $L^r$. 
Moreover, let $\nabla_f=\delta_f-\dbar$.

%%
%%if $\phi$ is a section of  $L^r$ and
%%$q=q_1+\ldots +q_m$ is a section of  $E\otimes L^r$, then
%%$\delta_f q=\phi$ means that $\sum_j p_jq_j =\phi$.

\begin{thm}\label{basic}
Let $f$ be holomorphic section of 
$E^*=L^{d_1}\oplus\cdots\oplus L^{d_m}$ and assume
that $\ell\ge 0$ is given and 
that 
$$
m-\ell\le n\quad {\rm or}\quad 
r\ge \sum_{j=1}^{n+\ell+1} d_j-n,
$$
where $d_1\ge d_2\ge\ldots\ge d_m$. 
If $\phi\in\O(\P^n,\Lambda^\ell E\otimes L^r)$,
then $\phi=\delta_f\psi$ for some
$\psi\in\O(\P^n,\Lambda^{\ell+1}E\otimes L^r)$
if and only if
\begin{equation}\label{gam}
\nabla_f (w\wedge R^f)=\phi\wedge R^f
\end{equation}
for some smooth $w$ defined in a   neighborhood of
$Z^f$.
\end{thm}

If $\ell>m-p$ ($p=\codim Z^f$)  then the condition on $\phi$  is void;
if $\ell=m-p$, it means that 
$\phi\wedge R^f=0$, see the remarks after Theorem~\ref{mfald} below.
If $f$ is a complete intersection, then 
$m\le n$ and therefore we have

\begin{cor} 
Let $f$ be a holomorphic section of 
$E^*=L^{d_1}\oplus\cdots \oplus L^{d_m}$ that   is a complete
intersection, and assume
that $r\ge 0$.
If $\phi\in\O(\P^n,L^r)$, then
$\phi=f\cdot q$ is solvable with
$q\in\O(\P^n,E\otimes L^r)$
if and only if
$\phi R^f=0$.
\end{cor}

\begin{proof}[Proof of Theorem~\ref{noll}]
If the hypotheses in Theorem~\ref{noll} are fulfilled,
then Theorem~\ref{basic} provides a section
$q=(q_1, \ldots , q_m)$ of  $E\otimes L^r$ such that
$\sum f_jq_j=\delta_fq=\phi$;
here $q_j$ are  sections  of $L^{-d_j+r}$.
After dehomogenization this means that 
$Q_j$ are polynomials %%with $\deg Q_j=r-d_j$
such that $\deg F_jQ_j\le r$.
\end{proof}

In Section~\ref{rest} we
recall the necessary background from \cite{A2} about the residue currents,
and present a general result about the image of a holomorphic
morphism $f$. Combined with well-known vanishing results
for the line bundles $L^r\to\P^n$ it leads to a proof of
Theorem~\ref{basic}.

In the last section we construct  explicit integral representations  of
the solutions in  Theorem~\ref{noll}.
They give essentially the same results except for  a small loss
of precision.
The construction is based on  ideas in \cite{A1} and \cite{A2}.

\section{The residue current of a holomorphic section}\label{rest}

Let $E\to X$ be a holomorphic Hermitian vector bundle of rank
$m$ over the $n$-dimensional complex manifold $X$,   and let $f$
be  a holomorphic section of  the dual bundle $E^*$, or
in other words,  a holomorphic morphism 
$f\colon E\to X\times\C$. 
Let 
$$
\L^r=\bigoplus_\ell \D'_{0,l+r}(X,\Lambda^\ell E);
$$
we consider $\L^r$ as a subbundle to $\Lambda(T^*_{0,1}\oplus E)$, so that
$\delta_f$ (i.e., interior multiplication with $f$) and $\dbar$ anticommutes.
Then  $\nabla_f=\delta_f-\dbar$  induces the complex
$\to\L^{r-1}\to\L^r\to$.
It is readily checked that $\nabla_f$ satisfies the Leibniz rule
$\nabla_f (\alpha\w\beta)=\nabla_f\alpha\w\beta+(-1)^\nu\alpha\w\nabla_f\beta$,
where $\nu$ is the total degree of $\beta$. 
%%%
Let $s$ be the dual section of  $E$ of $f$ so that in particular
$\delta_f s=\|f\|^2$. 
In \cite{A2} we defined  the current
$$
R^f=\dbar\|f\|^{2\lambda}\wedge \frac{s}{\nabla_f s}\big|_{\lambda=0};
$$
for large $\Re \lambda$ the right hand side is integrable and therefore 
a well defined
current, and by a nontrivial argument based on Hironaka's theorem
one can make an analytic continuation to $\lambda=0$. 
The resulting current  is an element in $\L^0$ with  support   on $Z^f=\{z;\ f(z)=0\}$
and it satisfies the basic equality
\begin{equation}\label{baseq}
\nabla_f U^f=1- R^f,
\end{equation}
where  $U^f\in\L^{-1}$ is defined as 
$$
U^f=\|f\|^{2\lambda} \frac{s}{\nabla_f s}\big|_{\lambda=0}.
$$
Moreover,
\begin{equation}\label{sesa}
R^f=R^f_{p,p}+\ldots + R^f_{m,m},
\end{equation}
where $p=\codim Z^f$; here lower index ${\ell,q}$ means that the current
has bidegree  $(0,q)$-form and takes values in
$\Lambda^\ell E$.

\begin{prop}\label{gubb}
Assume that  $f$ defines  a complete intersection and  that $h$ is a 
holomorphic section of some line bundle  such that $\{h=0\}\cap Z^f$ has codimension $m+1$.
If $\phi$ is a holomorphic section such that $\phi R^f=0$ in 
$X\setminus\{h=0\}$, then $\phi R^f=0$.
\end{prop}

Notice that since $f$ is a complete intersection, $R^f=R^f_m$.
The following lemma, which is the core of the proof, states that
then $R^f$ is  robust in a certain sense.

\begin{lma}\label{brott}
The current 
$|h|^{2\lambda}R^f$ has an  analytic continuation to $\Re \lambda>-\epsilon$ and
$$
|h|^{2\lambda}R^f|_{\lambda=0}=R^f.
$$
\end{lma}

\begin{proof}
Clearly the statement is local. By Hironaka's theorem and
a toric resolution we may assume that $f=f_0f'$,  where
$f_0$ is a holomorphic function and $f'$ is a non-vanishing
section.  In this way we can
 write the action of $R^f$ on a test form $\xi$
as a finite sum of terms like
$$
\int\dbar\Big[\frac{1}{f_0^\ell}\Big]\w\alpha\w\tilde\xi\rho,
$$
where $[1/f^\ell_0]$ is the principal value current,
$\alpha$ is a $(0,m-1)$-form, $\tilde\xi$ is the pull-back of
$\xi$ in the given resolution, and $\rho$ is a cut-off function. 
We may also assume that 
$$
f_0=\tau^{\alpha_1}_{k_1}\cdots\tau^{\alpha_\nu}_{k_\nu},
$$
in  appropriate local coordinates $\tau_j$, 
and therefore the integral is a sum of terms like
\begin{equation}\label{skrutt}
\int\Big[\Pi_{r\neq j}\frac{1}{\tau^{\alpha_r\ell}_{k_r}}\dbar\frac{1}{\tau^{\alpha_j\ell}_{k_j}}
\Big]\w \alpha\w\tilde\xi \rho.
\end{equation}
We may also assume that 
$
h=\tau_{m_1}^{\beta_1}\cdots\tau_{m_\mu}^{\beta_\mu} u,
$
where $u\neq 0$. 
Thus $|h|^{2\lambda}R.\xi$ is a finite sum of terms like
\begin{equation}\label{skrutt2}
\int|\tau_{m_1}|^{2\lambda\beta_1}\cdots|\tau_{m_\mu}|^{2\lambda\beta_\mu} |u|^{2\lambda}
\Big[\Pi_{r\neq j}\frac{1}{\tau^{\alpha_r\ell}_{k_r}}\dbar\frac{1}{\tau^{\alpha_j\ell}_{k_j}}
\Big]\w \alpha\w\tilde\xi \rho.
\end{equation}
If one of the $m_i$ is equal to $k_j$, then clearly this integral vanishes for
$\Re\lambda>>0$, and trivially therefore it has  an analytic continuation to
$\lambda>-\epsilon$, with the value $0$ at $\lambda=0$.  However, since
$\tau_{k_j}$ is a factor in both $h$ and $f_0$, and $\codim\{h=0\}\cap Z=m-1$,
for degree reasons it  follows that $\xi$ vanishes on this set, and therefore,
cf.\ e.g., \cite{BY2},  \cite{PTY} or \cite{A2}, each term  in $\tilde\xi$ contains either
a factor $\bar\tau_{k_j}$ or $d\bar\tau_{k_j}$. In any case, this implies that
already the integral \eqref{skrutt} vanishes. 
On the other hand, if no $m_i$ is equal to $\tau_{k_j}$,
it is easy to see that \eqref{skrutt2} has an analytic continuation to $\Re\lambda>-\epsilon$
and  takes   the  value   \eqref{skrutt}  at $\lambda=0$. %%% is equal to \eqref{skrutt}.
In fact, this follows easily  since if
$[1/s^\ell]$ is the usual principal value distribution in $\C$ and $v>0$ is smooth and strictly
positive, then
$$
|s|^{2\lambda} v^\lambda [1/s^\ell]
$$
has an analytic continuation to $\Re\lambda>-\epsilon$ and takes  the value 
$[1/s^\ell]$ at $\lambda=0$. Thus the proposition is proved.
\end{proof}

\begin{proof}[Proof of Proposition~\ref{gubb}]
By assumption $\phi R^f$ is a current with support on $\{h=0\}$, and hence (locally)
$|h|^{2\lambda}\phi R^f=0$ if $\Re\lambda>> 0$. From  Lemma~\ref{brott}  it follows that
$$
\phi R^f= |h|^{2\lambda}\phi R^f|_{\lambda=0}=0.
$$
\end{proof}

\smallskip

Let $L\to X$ be a holomorphic line bundle and let $\phi$ be a holomorphic section
of  $\Lambda^kE\otimes L$.

\begin{thm}\label{mfald}
Let $\ell\ge 0$ and
suppose that $H^{0,s}(X,\Lambda^{s+\ell+1}E\otimes L)=0$ for all 
$1\le  s \le m-\ell-1$. 
Moreover, let $\phi\in\O(X,\Lambda^\ell E\otimes L)$. Then
$\delta_f\psi=\phi$ has a  solution 
$\psi\in\O((X,\Lambda^{\ell+1}E\otimes L))$
if and only if there is a smooth solution $w$, defined
in a neighborhood of $Z^f$,   to
\begin{equation}\label{svada}
\nabla_f(w\wedge R^f)=\phi\wedge R^f.
\end{equation}
\end{thm}

In view of \eqref{sesa}, the condition on $\phi$  is void
if $\ell>m-p$.  Moreover, since 
$w=w_{\ell+1,0}+w_{\ell+2,1}+\cdots$ 
the condition means precisely that
$\phi\wedge R^f=0$ if $\ell=m-p$.
In the case $\ell=0$ and $p=m$, i.e., $f$ defines a complete
intersection, we get back the well-known duality theorem, first
proved in \cite{DS} and \cite{P1}.  %%In particular, when $c=m$, i.e., $f$ is a complete
%%intersection, then blablbla.
\smallskip

It was  also proved in \cite{A2} that $h^{\min(m,n)}R^f=0$ if $h$ is
holomorphic and $\|h\|\le C\|f\|$.
The local version of Theorem~\ref{mfald} therefore  immediately implies
 the Brian\c con-Skoda theorem, \cite{BS},:
\noindent {\it If $\|\phi\|\le C\|f\|$, then
locally $\phi^{\min(m,n)}$ belongs to the ideal $(f)$.}
There is also an explicit representation formula in \cite{A2}.

\begin{proof}[Proof of Theorem~\ref{mfald}]
First suppose that the holomorphic solution $\psi$ exists. Then
$\nabla_f\psi=\phi$ and hence
$\nabla_f(\psi\wedge R^f)= \phi\wedge R^f$
since $\nabla_f R^f=0$.
Conversely,
if \eqref{svada} holds for some smooth $w$, we claim that 
$\nabla_f v=\phi$, if 
$$
v= (-1)^{\ell}\phi\w U^f +w\wedge R^f.
$$
In fact, since
$\nabla_f\phi=0$, 
$$
\nabla v=\phi\w\nabla_f U^f+\nabla_f (w\w R^f)=
\phi\w (1-R^f)+ \phi\w  R^f=\phi.
$$
This means that 
%%if lower index denotes degree in $\Lambda^{\bullet}E$, then
%%this means that 
$$
\dbar v_{m,m-\ell-1}=0 \quad {\rm and}\quad \delta_f v_{k+1,k-\ell}=\dbar v_{k,k-\ell-1}.
$$
%%Notice that $v_{$ is a current of bidegree $(0,k-\ell-1)$.
By the assumption on the Dolbeault cohomology, we can successively solve
the equations
$$
\dbar \eta_{m,m-\ell-2}=v_{m,m-\ell-1},
\quad \dbar\eta_{k,k-\ell-2}=v_{k,k-\ell-1}+\delta_f\eta_{k+1,k-\ell-1},\quad k\ge \ell,
$$
and then finally  $\psi= v_{\ell,0}+\delta_f\eta_{\ell+1,0}$
is the  desired  holomorphic solution.
\end{proof}

\begin{ex}
Suppose that  $X$ is a compact 
and  $L$ is a strictly positive
line bundle. Then there is an $r_0>$ such that
$H^{0,k}(X,\Lambda E^{\bullet}\otimes L^r)=0$
for all $k\ge 1$ if $r\ge r_0$.
If $f$ is a holomorphic section of  $E^*$, 
then  a holomorphic section $\phi\in\O(\Lambda^\ell E\otimes L^r)$,
$r\ge r_0$, 
is in the image of the morphism
\begin{equation}\label{pil}
\O(X,\Lambda^{\ell+1}E\otimes L^r)\to \O(X,\Lambda^\ell E\otimes L^r)
\end{equation}
if $\phi\wedge R^f=0$. If $\ell=m-p$ the condition is necessary.
\end{ex}

\bigskip

We shall now focus on the case where  $X=\P^n$ and 
$E$ is the Hermitian vector bundle from Section~1. 
Let $E_1,\ldots, E_m$ be trivial line bundles over $\P^n$ with
basis elements $\epsilon_1,\ldots,\epsilon_m$, and let
$E_j^*$ be the dual bundles, with bases $\epsilon_j^*$.
Then we have that 
$$
E^*=(L^{d_1}\otimes E^*_1)\oplus\cdots  \oplus (L^{d_m}\otimes E^*_m),
$$
$$
E=(L^{-d_1}\otimes E_1)\oplus\cdots  \oplus (L^{-d_m}\otimes E_m),
$$
and  for instance our section $f$ can be written
$$
f=\sum_1^m f_j\epsilon^*_j.
$$
Its dual section $s$ is then, cf., \eqref{pnorm},
$$
s=\sum_j \frac{\overline{f_j(z)}}{|z|^{2d_j}}\epsilon_j,
$$
so
$$
R^f=\dbar\|f\|^{2\lambda}\wedge\sum_{\ell+1}^m\frac{s\wedge(\dbar s)^{\ell-1}}
{\|f\|^{2\ell}}\Big|_{\lambda=0}.
$$
In $\C^n=\{z_0\neq 0\}\subset\P^n$ we have the coordinates $z'$ and the natural
holomorphic frame $e_j=z_0^{-d_j}\epsilon_j$ and its dual
$e_j^*=z_0^{d_j}\epsilon_j^*$.
If $f'_j(z')=f_j(1,z')$ then
$$
f=\sum_1^m f'_j e_j^*
$$
and
$$
s=\sum_1^m \frac{\overline{f'_j(z')}}{(1+|z'|)^{d_j}} e_j.
$$
When $\codim Z^f=m$, the residue current $R^f$ is independent
of the metric, it  just contains the
top degree term $R^f_{m,m}$,  and in fact,    see \cite{A2}, 
$$
R^f=\big[\dbar\frac{1}{f'_m}\wedge\ldots\wedge\dbar\frac{1}{f'_1}\big]\wedge
e_1\wedge\ldots\w e_m,
$$
where the expression in brackets is a  Coleff-Herrera residue current.
Choosing the local coordinates $z_0, \zeta_1,\ldots,\zeta_n$ in
$\C^{n+1}\setminus\{0\}$, where   $\zeta_j=z_j/z_0$, it is easy to see that
$$
\pi^*\big[\dbar\frac{1}{f'_m}\wedge\ldots\wedge\dbar\frac{1}{f'_1}\big]=
z_0^{\sum d_j}\big[\dbar\frac{1}{f_m}\wedge\ldots\wedge\dbar\frac{1}{f_1}\big],
$$
%%(recall that the pull-back of a current is well-defined if the mapping
%%has surjective derivative)
and hence we can identify $R^f$ with the Coleff-Herrera current
$$ 
\big[\dbar\frac{1}{f_m}\w\ldots\w\dbar\frac{1}{f_1}\big]\w\epsilon_1\wedge\ldots
\w\epsilon_m
$$
in $\C^{n+1}\setminus\{0\}$.

%%
%%
%%If we  consider it as
%%a Coleff-Herrera current on $\C^{n+1}$, then it coincides 
%%on $\C^{n+1}\setminus\{0\}$ with the  pullback of $R^f$ to
%%$\C^{n+1}\setminus\{0\}$ under the natural projection $\C^{n+1}\setminus\{0\}\to
%%\P^n$, i.e., 

\begin{proof}[Proof of Theorem~\ref{basic}]
It is well-known, see, e.g., \cite{Dem}, that
$H^{0,k}(\P^n, L^\nu)=0$ for all $\nu$ if 
$1\le k\le n-1$ and that
$H^{0,n}(\P^n, L^\nu)=0$ if (and only if) $\nu\ge -n$.
Since   $E=L^{-d_1}\oplus\cdots\oplus L^{-d_m}$ we have that 
$$
\Lambda^\nu E\otimes L^r=\bigoplus'_{|J|=\nu}L^{-d_{J_1}}\otimes\cdots
\otimes L^{-d_{J_{\nu}}}\otimes L^r=
\bigoplus'_{|J|=\nu}L^{r-d_{J_1}\cdots-d_{J_\nu}}.
$$
Thus $H^{0,s}(\P^n,\Lambda^{s+\ell+1}E\otimes L^r)=0$ for 
$1\le s\le m-\ell-1$ if  either $m-\ell-1\le n -1$ or 
$$
r-\sum_1^{n+\ell+1} d_j \ge -n.
$$
Now Theorem~\ref{basic} follows from
Theorem~\ref{mfald}. %%% with $f=p$.
\end{proof}

\section{Integral representation}

The aim of this section is to present an explicit integral representation
of the solution $Q_j$ to the division problem in Theorem~\ref{noll}.
We have

\begin{thm}\label{skrott}
 Let $F_1,\ldots, F_m,\Phi$ be  polynomials  in $\C^n$, 
let $f$ and  $R^f$ be as before, and let $\phi$  be the $r$-homogenization
of $\Phi$ ($\deg\Phi\le r$).  Then there is an explicit
decomposition
\begin{equation}\label{deco}
\Phi(z')=\sum_1^m F_j(z')\int_{\P^n} T^j(\zeta,z')\phi(\zeta)
+\int_{\P^n} S(\zeta,z')\wedge R^f(\zeta)\phi(\zeta),
\end{equation}
where $T^j(\zeta,z'), S(\zeta,z')$ are smooth forms (in $[\zeta]$) on $\P^n$ and  
holomorphic polynomials in $z'$,   such that 
$$
\deg_{z'}\big(F_j(z')T^j(\zeta,z'))\le 
d_1+d_2+\cdots +d_{\mu+1}+r,
$$
if $\mu=\min(n,m-1)$ and $d_1\ge d_2\ge\cdots\ge  d_m$.
\end{thm}

Thus, if $\phi R^f=0$ we get back the conclusion of Theorem~\ref{noll} but with the
extra term $d_1+\cdots+d_{\mu+1}$ in the estimate of the degree.

%%Our starting point is the construction of division formulas for trivial
%%bundles over domains in $\C^n$ in \cite{A2}, 
%%which  we adapt to the present case.

\smallskip
For fixed $z\in\C^n$, 
$$
\eta=2\pi i\sum_0^n z_j \frac{\partial}{\partial\zeta_j}
$$ 
is an  
$L_z\otimes L^{-1}_\zeta$-valued 
$(1,0)$-form on $\P^n$,
and if $\delta_\eta$ denotes interior multiplication with $\eta$, then
$$
\delta_\eta\colon \D'_{\ell+1,0}(\P^n,  L^{r+1})\to 
\D'_{\ell,0}(\P^n, L^{r}).
$$

\begin{remark}
When we say that $\eta$ is a section of  $L_z\otimes L^{-1}_\zeta$ rather than
 $L^{-1}=L^{-1}_\zeta$, we just indicate that it is $1$-homogeneous
in $z$; it would be more correct, but  less convenient,
to consider  $\eta$ as   a section of  the bundle 
$L_z\otimes L^{-1}_\zeta\otimes (T^*_\zeta)_{0,1}$ over
$\P^n_z\times\P^n_\zeta$.
\end{remark}

Let $\nabla_\eta=\delta_\eta-\dbar$.
Notice that if 
$$
\alpha=\alpha_0+\alpha_1= \frac{z\cdot\bar \zeta}{|\zeta|^2}-
\dbar\frac{\bar\zeta\cdot d\zeta}{2\pi i|\zeta|^2},
$$
then the first term, $\alpha_0$,  is a section of  $L_z\otimes L_\zeta^{-1}$ and
the second term,  $\alpha_1$,  is a projective form (since
$\delta_\zeta \alpha_1=0$); moreover 
\begin{equation}\label{alphaeq}
\nabla_\eta \alpha=0.
\end{equation}
We have the following basic integral representation of global holomorphic sections
of  $L^r$.
\begin{prop}\label{grodor}
Assume that $r\ge 0$ and that $\phi\in\O(\P^n, L^r)$. Then
$$
\phi(z)=\int_{\P^n}\alpha^{n+r}\phi.
$$
\end{prop}

For degree reasons, actually
$$
\phi(z)=\frac{(n+r)!}{n!r!}\int_{\P^n}\alpha_0^r\w \alpha_1^n\phi;
$$
this formula appeared already in \cite{BBB};
expressed in affine coordinates it is the well-known
weighted Bergman representation formula for polynomials
in $\C^n$. However, we prefer to  supply a  direct proof
on $\P^n$, following the ideas in \cite{A1}.

\begin{proof}
%%Expressed in affine coordinates this is a well-known representation formula
%%for polynomials in $\C^n$, but we prefer to 
Let $\sigma$ be the  $L^{-1}_z\otimes L_\zeta\otimes T^*_{1,0}(\P^n_\zeta)$
valued $(1,0)$-form on $\P^n$ 
that is dual, with respect to the natural metric, to $\eta$.
Then, since $\eta$ has a first order zero  at $[z]$ (and no others), it follows 
(see \cite{A1})   that
$$
\nabla_\eta\frac{\sigma}{\nabla_\eta\sigma}=1-\big[[z]\big].
$$
The rightmost term  is the  $L_z^{-n}\otimes L^n_\zeta$-valued   $(n,n)$-current 
point evaluation at $[z]$ for  sections of  $L^{-n}$.
If $\phi$ is a global holomorphic section of  $L^r$ it follows by \eqref{alphaeq} that
$$
\nabla_\eta\big(\frac{\sigma}{\nabla_\eta\sigma}\wedge \alpha^{n+r}\phi\big)=
\phi\alpha^{n+r} -\phi\big[[z]\big],
$$ 
where this time the last term is $\phi$ times the $L^r_z\otimes L^{-r}_\zeta$-valued
current point evaluation at $[z]$.
If we integrate this equality over $\P^n$ we get the desired
representation formula.
\end{proof}

Let $E_1,\ldots, E_m$ be the trivial line bundles over $\P^n$ with
basis elements $\epsilon_1,\ldots,\epsilon_m$, 
so that $E=(L^{-d_1}\otimes E_1)\oplus\cdots  \oplus (L^{-d_m}\otimes E_m)$
as in Section~\ref{rest}.
We also introduce   disjoint copies $\widetilde E_j$ of $E_j$ with bases
$\tilde \epsilon_j$ and the bundle
$$
\widetilde E=
(L^{-d_1}\otimes \widetilde E_1)\oplus\cdots  \oplus (L^{-d_m}\otimes \widetilde E_m).
$$
Let  $\Lambda$ be the exterior algebra bundle over 
the direct sum of all the bundles
$E$, $\widetilde E$, $E^*$, and  $T^*(\P^n)$.
Any form $\gamma$ with values in $\Lambda$    can be written
uniquely as $\gamma=\gamma'\wedge (\sum \epsilon_j^*\w \epsilon_j)^m/m! + \gamma''$
where $\gamma''$ denotes terms that do not contain a factor 
$(\sum \epsilon_j^*\w \epsilon_j)^m/m!$, and we define
$$
\int_\epsilon \gamma=\gamma'.
$$
We have a globally defined form 
$$
\tau=\sum_1^m \epsilon_j^*\wedge(\epsilon_j-\tilde \epsilon_j).
$$
%%and if $\gamma$ is a section of  $\Lambda E$
%%and $\tilde\gamma=\sum\gamma_j\tilde\epsilon_j$ 
%%is the  corresponding section of  $\widetilde E$,
%%then,  see \cite{A2},  
%%\begin{equation}\label{sprut}
%%\int_\epsilon \tau_m\wedge \gamma=\tilde \gamma.
%%\end{equation}

\smallskip

From now on we consider $[z]$ as a fixed arbitrary point in $\C^n\subset\P^n$, and
let $z=(1,z')$.
%%
%%and define the two sections 
%%$$
%%f=\sum_j z_0^{d_j}p_j(\zeta)\epsilon_j
%%$$
We also  introduce the section 
$$
f_z=\sum_j \zeta_0^{d_j}f_j(1,z)\epsilon^*_j=\sum\zeta_0^{d_j}F_j(z')\epsilon^*_j
$$
of  $E^*$ and let $\tilde f_z$ be the corresponding section of  $\tilde E^*$.

%%\begin{remark}
%%It might  have been more natural to consider $E^*$ as the bundle
%%$$
%%E^*=L_z^{d_1}\otimes L^{d_1}\otimes E^*_1\oplus 
%%=L_z^{d_2}\otimes L^{d_2}\otimes E^*_2\cdots
%%$$
%%over $\P^n\zeta\times\P^n_z$
%%but for technical reasons we keep   $z$ just as a parameter.
%%\end{remark}

\begin{lma}
There is a holomorphic section $H=\sum H_j\w \epsilon_j$
of  $E^*\otimes L\otimes T^*_{1,0}$, thus  
$H_j$ are sections of  $L^{d_j}\otimes L\otimes T^*_{1,0}$,
such that
$$
\delta_\eta H= f- f_z,
$$
and such that the coefficients in $H_j$ are
polynomials in $z'/z_0$ of degrees (at most)
$d_j-1$.
\end{lma}

\begin{proof}
For each $F_j(z')$ we can find Hefer functions $h_j^k(\zeta',z')$,
polynomials of degree $d_j-1$ in $(\zeta', z')$, such that
$$
\sum_{k=1}^nh_j^k(\zeta',z')(\zeta_k-z_k)=F_j(\zeta')-F_j(z').
$$
If we then take
$$
H_j=\frac{\zeta_0^{d_j+1}}{2\pi i}
\sum_1^n h_j^k(\zeta'/\zeta_0, z')d(\zeta_k/\zeta_0),
$$
then clearly $H_j$ is a  projective $(1,0)$-form, and moreover,
$$
\delta_\eta H_j=f_j(\zeta)-\zeta_0^{d_j}F_j(z')
$$
as wanted.
\end{proof}

Let  $\delta_F$ %%=\delta_f+\delta_{\tilde{\check f}}$ d
denote  interior multiplication with the section
 $F=f+\tilde f_z$ of   $E^*\oplus \tilde E^*$. Then
$\delta_F \tau=  f-f_z= -\delta_\eta H$.
If 
$$
\nabla=\delta_F+\delta_\eta-\dbar,
$$
thus
\begin{equation}\label{ht}
\nabla(\tau+H)=0.
\end{equation}
We are now ready to define the explicit division formula.

%%\begin{prop}\label{gropp}
%%Let $\phi$ be an $L^r$-valued holomorphic section.
%%With the notation above we have that
%%
%%\begin{multline}\label{skot}
%%\phi(1,z')=
%%\int_{\P^n}\int_\epsilon \delta_{\tilde p_z}\big[e^{\tau+H}\wedge  
%%U^f\wedge \alpha^{n+r}\wedge \phi\big] +\\
%% + \int_e\int_X e^{\tau+H}\wedge  U\wedge g \wedge \delta_f\phi 
%%+ \int_{\P^n}\int_\epsilon e^{\tau+H}\wedge R^f\wedge \alpha^{n+r} \wedge \phi.
%%\end{multline}
%%\end{prop}

\begin{proof}[Proof of Theorem~\ref{skrott}]
From \eqref{ht} it follows that
\begin{equation}
(\nabla_{\eta}+\delta_F)(e^{\tau+H}\wedge U^f)= e^{\tau+H}\wedge(1-R^f).
\end{equation}
We can rewrite this as
\begin{equation}\label{snabel}
\delta_F(e^{\tau+H}\wedge U^f)+e^{\tau+H}\wedge R^f=
e^{\tau+H}-\nabla_{\eta}(e^{\tau+H}\wedge U^f).
\end{equation}
We claim that the component of full bidegree $(n,n)$ of 
\begin{equation}\label{snutte}
\int_{\epsilon}
\big[e^{\tau+H}-\nabla_{\eta}(e^{\tau+H}\wedge U^f)\big]\wedge \alpha^{n+r}\phi
\end{equation}
is equal to 
$$
\frac{(n+r)!}{n!r!}\alpha_1^n\alpha_0^r\phi+\dbar(\cdots)
$$
where $(\cdots)$ is a scalar-valued $(n,n-1)$-form.
In fact, since $\alpha^{n+r}$ has bidegree $(*,*)$ the factor
$U_{\ell,\ell-1}$ must be combined with $H_\ell$, and then it follows that
$\tau$ can be replaced by $\omega=\sum_j\epsilon^*_j\w\epsilon_j$.
Observe that the component of $U_{\ell,\ell-1}$ with basis element
$\epsilon_{J_1}\w\ldots\w\epsilon_{J_\ell}$ takes values in
$L^{-(d_{J_1}+\cdots +d_{J_{\ell}})}$, whereas the component of
$H_\ell$ with basis element $\epsilon^*_{J_1}\w\ldots\w\epsilon^*_{J_\ell}$
takes values in $L^{d_{J_1}+\cdots +d_{J_{\ell}}}\otimes L^\ell$.
The product of these two factors must be combined with 
$\alpha_1^{n-\ell}\alpha_0^{\ell+r}\phi$ which gives 
a scalar-valued $(n,n)$-form as claimed. 
Thus we can integrate \eqref{snutte}  over $\P^n$, 
and by Proposition~\ref{grodor} and  Stokes' theorem 
it is equal to   $\phi(z)$.

We now consider the left hand side of \eqref{snabel}
multiplied with $\alpha^{n+r}\phi$.
To begin with,
$$
\int_{\P^n}\int_\epsilon e^{\tau+H}\wedge R^f
\wedge \alpha^{n+r}\phi
$$
is well defined with the same argument as above, and again
one can replace $\tau$ by $\omega$.
Moreover, since $\alpha^{n+r}\phi$ contains no $\epsilon_j$,
$$
\int_\epsilon\delta_f(e^{\tau+H}\wedge U^f)
\wedge \alpha^{n+r}\phi=
\int_\epsilon\delta_f(e^{\tau+H}\wedge U^f
\wedge \alpha^{n+r}\phi)=0.
$$
Since
$$
\delta_{\tilde f_z}\sum_j\tilde\epsilon_j\wedge\epsilon^*_j
=\sum_j F(z')\zeta_0^{d_j}\epsilon_j^*=f_z,
$$
another computation shows that the component of bidegree $(n,n)$ of
$$
\int_\epsilon \delta_{\tilde f_z}(e^{\tau+H}\wedge U^f)
\wedge \alpha^{n+r}\phi
$$
is equal to
$$
\int_\epsilon
f_z\wedge\sum_{k=0}^{m-1}\omega_{m-k-1}\wedge H_k\wedge U_{k+1,k}\wedge\alpha_1^{n-k}
\alpha_0^{k+r}\phi.
$$
Again one can check that this form is scalar valued.
Summing up we have the desired  decomposition \eqref{deco}
with
\begin{multline*}
S(\zeta,z')\wedge R^f(\zeta)=
\int_\epsilon  e^{\omega+H}\wedge R^f
\wedge \alpha^{n+r}=\\
\sum_{k=\codim Z^f}^m\int_\epsilon \frac{(n+r)!}{(n-k)!(k+r)!}
\omega_{m-k}\w H_k\w R^f_{k,k}\alpha_1^{n-k}\alpha_0^{k+r},
\end{multline*}
and
\begin{multline*}
T^j(\zeta,z')=\\
\int_\epsilon
\epsilon_j^*\zeta_0^{d_j}
\wedge\sum_{k=1}^{m-1}
\frac{(n+r)!}{(n-k)!(k+r)!}
\tilde I_{n-k-1}\wedge H_k\wedge U_{k+1,k}\wedge
\alpha_1^{n-k}\alpha_0^{k+r}\phi,
\end{multline*}
Both  $\alpha$ and $H$ are polynomials in $z'$
so it just remains to check the  degrees of $T^j$. 
The worst case occur when $k$ is as large as possible which is
$k=\mu=\min(m-1,n)$. Then the factor $\alpha_0^{k+r}$
has degree $k+r$. 
Recall that  $H=\sum H_\ell\w\epsilon_\ell^*$ 
and that $\deg H_\ell=d_\ell-1$. The term  $H_j$ cannot occur,
because of the presence of $\epsilon^*_j$, 
and thus we get that $d_j+\deg Q_j$ is at most 
$d_1-1 +d_2-1 +\cdots d_{\mu+1} -1 +1 +\mu+r=d_1+\cdots +d_{\mu+1}+r$. 
\end{proof}

The division formula  constructed here, Theorem~\ref{skrott},  is a generalization to $\P^n$ of the formula in
\cite{A2},  which was used to give an explicit representation of  the solutions  in the local
version of Theorem~\ref{mfald}; in particular it provided the first known explicit
proof  of the Brian\c con-Skoda theorem.
This division formula is  based on the ideas in \cite{A1} and it
differs from Berndtsson's  classical formula, \cite{BB}, 
in some  respects. To begin with our formula  works also for sections with values in
$\Lambda^\ell E$, although in this paper we have only generalized
the scalar-valued part  to $\P^n$.
The more interesting novelty %%of the formula in \cite{A2},  
with regard to this paper,
is that the residue term contains precisely the factor $\phi R^f$, so that our   formula 
provides a solution of the division problem as soon as $\phi R^f=0$ (or $\phi R^f=\nabla_f(w\w R^f)$
for some smooth $w$). 
One can obtain  a similar formula involving residues 
(but not precisely $R^f$ except for the complete intersection case)
from  Berndtsson's formula;  this was first
done by Passare in \cite{P1}, and  various variants
have  been used by several authors since then, see
\cite{BGVY} and the references given there. 
These formulas all go back to the construction of weighted
integral formulas in \cite{BA}.
However, the division formula in \cite{A2}, even in the simplest case,
when $f$ is nonvanishing,  could not  have  been obtained from
\cite{BA}, because the required choice of weight,
see formula (2.12) in Remark~3 in \cite{A3}, is not 
encompassed by the method in \cite{BA}, but the more general construction in
\cite{A1} is needed.

\smallskip

\noindent{\bf Ackowledgement}
I am indebted to the referee for his careful reading and for his  many important remarks
and constructive suggestions that have helped to clarify and improve the  final version of
this  paper. %manuscript.

\def\listing#1#2#3{{\sc #1}:\ {\it #2},\ #3.}

\end{document}